\documentclass[12pt]{amsart}
\usepackage{amsmath}
\usepackage{amsfonts,amssymb,amsthm, txfonts, pxfonts}
\newcommand{\R}{{\bf R}}

\newcommand{\g}{{\mathfrak g}}

\newcommand{\h}{{\mathfrak h}}

\newcommand{\ad}{\text{ad}}
\newcommand{\Ad}{\text{Ad}}

\newcommand{\Trace}{\text{Trace}}
\newtheorem{theorem}{{\sc Theorem}}
\newtheorem{cor}{{\sc Corollary}}
\newtheorem{lemma}{{\sc Lemma}}
\newcommand{\rf}[1]{(\ref{#1})}

\newcommand{\Sym}{\mathcal{S}}
\newcommand{\F}{\mathcal{F}}
\newcommand{\D}{\mathcal{D}}
\newcommand{\E}{\mathcal{E}}
\renewcommand{\L}{\mathcal{L}}

\newcommand{\beq}{\begin{equation}}
\newcommand{\eeq}[1]{\label{#1}\end{equation}}
\newcommand{\bea}[1]{\begin{array}{#1}}
\newcommand{\ea}{\end{array}}
\newcommand{\inter}{\cap}

\begin{document}

\title{Differentiability of Functions of Matrices}
\shorttitle{Matrix Functions}
\author{Yury Grabovsky, Omar Hijab \& Igor Rivin}
\address{Department of Mathematics, Temple University, Philadelphia, PA 19122}
\email{yury@math.temple.edu}
\email{ hijab@math.temple.edu}
\email{ rivin@math.temple.edu}

\thanks{Yury Grabovsky gratefully acknowledges the
support of the National Science Foundation through the grants NSF-0094089
and NSF-0138991.  Igor Rivin gratefully acknowledges the support of
the National Science Foundation through DMS-0072622}

\date{October 6, 2003}


\begin{abstract}
Let $f$ be a function on the set of diagonal $n\times n$ matrices, and
let $\tilde{f}$ be the unique extension of $f$ to the set of symmetric
$n\times n$ matrices invariant with respect to conjugation by
orthogonal matrices. We show that $\tilde{f}$ has the same regularity
properties as  $f.$ That is, if $f$ is $C^k,$ or $C^{k+\alpha},$ or
$C^\infty$ or $C^\omega$ than so is $\tilde{f}.$
\end{abstract}

\maketitle

It is well-known that every  rotation-invariant function $F$ on the space  $\Sym$ of real $d\times d$ symmetric matrices
 is determined by its restriction $f$ to the diagonal matrices,
$$f(r_1,\dots,r_d)=
F\begin{pmatrix}
r_1&0&0&\dots&0\\
0&r_2&0&\dots&0\\
&&\ddots&&\\
0&0&0&\dots&r_d
\end{pmatrix}.
$$
Since a $90^\circ$ rotation in the $ij$-plane interchanges $r_i$ and $r_j$, $f$ must necessarily be symmetric,
$$f(r_{\sigma1},r_{\sigma2},\dots,r_{\sigma d})=f(r_1,\dots,r_d),\qquad \text{for all permutations }\sigma.$$
It is then natural to seek properties of $f$ that are inherited by $F$. 

For example, suppose $f$ is a polynomial; then \cite{MR2003e:00003} 
$f=p(n_1,\dots,n_d)$ for some other polynomial $p$, where
$$n_k(r_1,\dots,r_d)=r_1^k+\dots+r_d^k,\qquad k\ge1,$$
are the Newton sums. It follows that
$$F(x)=p(\Trace(x),\dots,\Trace(x^d)),\qquad x\in \Sym,$$
since both sides are rotation-invariant and they agree on the diagonal matrices. Thus $f$ polynomial
implies $F$ polynomial.

Another interesting property is differentiability.
Recently, Lewis and Sendov \cite{MR2002i:15013} showed that if $f$ is $C^1$ or $C^2$, then  $F$ is 
$C^1$ or $C^2$ respectively; moreover they derived formulas for $DF(x)$ and $D^2F(x)$ in terms
of spectral quantities, i.e. the eigenvalues and eigenprojections of $x$. 
In this paper, we extend this result to $C^n$ and derive a formula  for $D^nF(x)$ in terms of the spectral 
quantities of $x$.

A theorem of C. Davis \cite{MR19:832b} asserts that $f$ convex implies $F$ convex, the canonical example
being the negative of the
logarithm of the determinant. There are several alternate proofs of this result, by Lewis \cite{MR2001g:49026}, 
Rivin \cite{math.FA/0208223}, 
and Grabovsky and Hijab \cite{GH}. As noted in \cite{MR2002i:15013}, this convexity result, in the 
$C^2$ setting, is
a consequence of the above differentiability result and the characterization of convexity in terms of 
nonnegativity of the second derivative.

These questions have natural generalizations in the context of compact Lie algebras.
In this setting the issue is to identify the interesting properties that are inherited by an
$\Ad$-invariant function $F$ on a compact Lie algebra $\g$ 
from its restriction $f$ to a Cartan subalgebra $\h$.  The polynomial question in this
setting is a theorem of Chevalley \cite{MR2003c:22001}, and the convexity question was extended
to this setting by Lewis \cite{MR2001g:49026} and subsequently by Grabovsky and Hijab \cite{GH}. 

For motivation, in \S1  we derive the analog of this result in the radial setting
and in \S2 we derive the result in the context of symmetric matrices.

\begin{section}{The Radial Case}

If $F:\R^d\to\R$ is continuous, then 
\beq
f(r,\pi)=F(x)\,\qquad x=r\pi,
\eeq{Ff}
is  continuous on $\R\times{\bf S}^{d-1}$, even, and $f(0,\pi)$ does not depend on $\pi$. 
Conversely, if $f:\R\times{\bf S}^{d-1}\to\R$ is  continuous, even, and $f(0,\pi)$ does not depend on $\pi$,
then
\beq
F(x)=f\left(|x|,\frac{x}{|x|}\right),\qquad x\not=0,
\eeq{fF}
extends to a continuous function on $\R^d$ satisfying \rf{Ff}. 

Let $\delta_\xi(\pi)=\xi-\langle \pi,\xi\rangle\pi$; 
then, for each $\xi$, the map $\delta_\xi:\R^d\to\R^d$ is a
vector field tangent to ${\bf S}^{d-1}$. If $f$ extends to a function on $\R\times\R^d$ that is polynomial
in $\pi$, then so does $\delta_\xi(f)$. For $f$ continuous in $r$ and polynomial in $\pi$, define 
$$\L_\xi^j(f)(r,\pi)=\langle\pi,\xi\rangle f(r,\pi)+\int_0^1 t^j \delta_\xi(f)(tr,\pi)\,dt.$$
Let $f'$ denote the derivative with respect to $r$.

For $x\not=0$, the maps $x\mapsto r=|x|$ and $x\mapsto \pi=x/r$ are analytic and their derivatives
in the $\xi$ direction are
$$r_\xi=\langle\pi,\xi\rangle,\qquad \pi_\xi=\frac{\delta_\xi(\pi)}{r}.$$
By \rf{Ff}, $\delta_\xi(f)(0,\pi)=0$ since $f(0,\pi)$ does not depend on $\pi$; 
if $f$  is $C^1$ in $r$ and polynomial in $\pi$, \rf{fF} and the chain rule implies
\beq
\begin{split}
D_\xi F(x)&=\langle\pi,\xi\rangle f'(r,\pi)+\frac{\delta_\xi(f)(r,\pi)}{r}\\
&=\langle\pi,\xi\rangle f'(r,\pi)+\int_0^1 \delta_\xi (f')(tr,\pi)\,dt=\L^0_\xi(f')(r,\pi),\quad x=r\pi\not=0.
\end{split}
\eeq{DF}
If $F$ is $C^1$ on $\R^d$, \rf{DF} is valid on $\R^d$; if $f$ is $C^2$ in $r$ and polynomial in $\pi$, 
we may repeat
this argument with $D_\xi F$ replacing $F$ and $\L^0_\xi(f')$ replacing $f$;
we obtain $D_\xi^2F(x)=\L_\xi^0(\L_\xi^0(f'))'(r,\pi)$ on $x=r\pi\not=0$. If $F$ is $C^{n-1}$ on $\R^d$ 
and $f$ is $C^n$ in $r$, we may continue in this manner to obtain
\beq
D_\xi^nF(x)=\L^0_\xi\L^1_\xi\dots\L^n_\xi\left(f^{(n)}\right)(r,\pi),\qquad x=r\pi\not=0;
\eeq{L}
here we used $(\L_\xi^jf)'=\L_\xi^{j+1}(f')$.

Now suppose $F$ is rotation-invariant; then $f=f(r)$ does not depend on $\pi$ hence \rf{L} implies
\beq
|D^n F(x)|\le C\sup_{|r|\le|x|}|f^{(n)}(r)|,\qquad |x|\not=0.
\eeq{ineq}

\begin{theorem}
Let $n\ge0$ and let $F$ be a rotation-invariant function on $\R^d$.
If the restriction  of $F$ to an axis is $C^n$, then $F$ is $C^n$ on $\R^d$.
\end{theorem}

The proof here mimics that of the matrix case in the next section; a simpler proof is possible.
\begin{proof}
Since $F$ is continuous on $\R^d$, we may derive this by induction, so we may assume $F$ is $C^{n-1}$
on $\R^d$. Since $f$ is even, $p(r)=f^{(n)}(0)r^n/n!$ either vanishes or is an even-order polynomial; 
hence $P(x)=p(|x|)$ is a polynomial on $\R^d$. Replacing $F$ by $F-P$,
we may further assume $f^{(n)}(0)=0$. In this case, by \rf{ineq}, we conclude $D^nF(x)\to0$ as $|x|\to0$. 
Since $F$ is $C^{n-1}$ on $\R^d$ and $C^n$ away from the origin, this implies $F$ is $C^n$ on $\R^d$.
\end{proof}

\end{section}

\begin{section}{The  Matrix Case}

Let $\Sym$ denote the vector space of real $d\times d$ symmetric matrices and let $G$ be the group
of $d\times d$ rotation matrices.   A function $F:\Sym\to\R$ is {\it rotation-invariant} if 
$F(gxg^{-1})=F(x)$ for every $x\in\Sym$ and $g\in G$. The result is

\begin{theorem}
\label{diffth}
Let $n\ge0$ and let $F$ be a rotation-invariant function on $\Sym$. 
If the restriction  of $F$ to the diagonal matrices $\D$ is $C^n$, then $F$ is $C^n$ on $\Sym$.
\end{theorem}

At the end of this section, we exhibit a formula \rf{der}  expressing $D^n F$ in terms of
derivatives of the restriction.

\begin{cor}
Let $F$ be a rotation-invariant function on $\Sym$ and let $f$ be its restriction to $\D$.
If $f$  is $C^\infty$ on $\D$, then $F$ is $C^\infty$ on $\Sym$.
If $f$  is $C^{n,\alpha}$, $0<\alpha<1$, on $\D$,
then $F$ is $C^{n,\alpha}$ on $\Sym$. If $f$  is analytic on $\D$, then $F$ is analytic on $\Sym$.
\end{cor}

The proof is at the end of the section.

A function $F$ on a vector space sum $A\oplus B$ is $C^{n,N}$ if the partial derivatives 
$D^{\alpha}_aD^\beta_bF$ exist and are continuous on $A\oplus B$ 
for all multi-indices $|\alpha|\le n$, $|\beta|\le N$.
The previous theorem is a special case of the slightly stronger 

\begin{theorem}Let $\E$ be a euclidean space and let $F:\Sym\oplus \E\to\R$ be rotation-invariant in the 
first variable. If the restriction of $F$ to $\D\oplus \E$ is $C^{n,N}$, then $F$ is $C^{n,N}$ on $\Sym\oplus \E$.
\end{theorem}

Let $\Sym_0$ denote the traceless matrices in $\Sym$, and let $\D_0=\D\inter\Sym_0$. 
Since $\Sym=\Sym_0\oplus\R$ and the trace is unchanged under conjugation by elements of $G$,  this in turn follows from

\begin{theorem}
\label{main}
Let $\E$ be a euclidean space and let $F:\Sym_0\oplus \E\to\R$ be rotation-invariant in the
first variable. If the restriction of $F$ to $\D_0\oplus \E$ is $C^{n,N}$, then $F$ is $C^{n,N}$ on $\Sym_0\oplus \E$.
\end{theorem}

We turn to the proof of Theorem \ref{main}.
To simplify notation, we now drop the subscript $0$, i.e. henceforth 
the spaces of traceless symmetric and diagonal matrices will be denoted
$\Sym$ and $\D$ respectively.
We will argue by induction over $(d,n,N)$, where we impose the lexicographic ordering on triples $(d,n,N)$.
{\it 
Thus we assume the result is true for all dimensions lower than $d$ and all orders of differentiability
in $n$ and $N$,
and we assume the result is true for dimension $d$ and all orders of differentiability on $\E$
and lower than $n$ on $\Sym$.  }
Let $\Sym^*$ be the open set of nonzero traceless symmetric matrices. 

Let $x_0\in\Sym^*$. Since conjugation by a rotation is
an invertible analytic map on $\Sym$, we may assume that $x_0$ is diagonal with the diagonal entries of $x_0$ arranged 
in decreasing order. Let $r_i(x_0)$, $i=1,\dots,d$, be the diagonal entries of $x_0$. 

Let $\g$ be the vector space of real $d\times d$ skew-symmetric matrices; then $\Sym\oplus\g$ is 
the vector space of all real $d\times d$ traceless matrices. 
Let $(\Sym\oplus\g)_0$ be the subspace of matrices $x$ commuting with $x_0$,
$xx_0=x_0x$. Then $x\in (\Sym\oplus\g)_0$ iff $x_{ij}=0$ whenever $r_i(x_0)\not=r_j(x_0)$. Let $G_0=G\inter
(\Sym\oplus\g)_0$, $\Sym_0=\Sym\inter(\Sym\oplus\g)_0$, and $\g_0=\g\inter(\Sym\oplus\g)_0$. 
Then $G_0$ is the isotropy group  of $x_0$ under the conjugation action, $(\Sym\oplus\g)_0=\Sym_0\oplus\g_0$,
and matrices in $\Sym_0$, $\g_0$, and $G_0$ are block-diagonal with 
the same block structure. If $[x,y]=xy-yx$ is the usual bracket, then
$\g_0=[\Sym_0,\Sym_0]$ and $\Sym_0$ is the orbit of $\D$ under conjugation by matrices in $G_0$.

\begin{lemma}
\label{S0}$F$ is $C^{n,N}$ on $\Sym_0\oplus\E$.
\end{lemma}

\begin{proof}
We derive this by applying the inductive hypothesis block-by-block.
Let $\Sym_0^{(k)}$ be the vector space of traceless symmetric matrices  in $\Sym_0$ where all blocks, 
except possibly the $k$-th block, vanish, and let $\D^{(k)}=\Sym_0^{(k)}\inter \D$.
Then $\Sym_0^{(1)}\oplus(\D\ominus\D^{(1)})\oplus\E$ consists of block-diagonal matrices that are
diagonal in all but the first block. Since the dimensions of the first block are strictly less than
$d$ and $F$ is $C^{n,N}$ on $\D$, the inductive hypothesis implies $F$ is $C^{n,N}$ on 
$\Sym_0^{(1)}\oplus(\D\ominus\D^{(1)})\oplus\E$. 
More precisely, if $\E^{(1)}=(\D\ominus\D^{(1)})\oplus \E$, then $F$  is $C^{n,N}$ on 
$\D\oplus\E=\D^{(1)}\oplus\E^{(1)}$. Since the matrices in $\D^{(1)}$ are strictly smaller than $d\times d$,
by the inductive hypothesis, $F$ is $C^{n,N}$ on
$$\Sym_0^{(1)}\oplus \E^{(1)}=\Sym_0^{(1)}\oplus(\D\ominus\D^{(1)})\oplus\E.$$
If $\E^{(2)}= \Sym_0^{(1)}\oplus(\D\ominus(\D^{(1)}\oplus\D^{(2)}))\oplus\E$, decomposing 
$\Sym_0^{(1)}\oplus(\D\ominus\D^{(1)})\oplus\E$ into $\D^{(2)}\oplus\E^{(2)}$ and applying 
the inductive hypothesis again, $F$ is $C^{n,N}$ on 
$$\Sym_0^{(2)}\oplus\E^{(2)}
=(\Sym_0^{(1)}\oplus\Sym_0^{(2)})\oplus(\D\ominus(\D^{(1)}\oplus\D^{(2)}))\oplus\E.$$
Continuing in this manner, we conclude $F$ is $C^{n,N}$ on 
$$\Sym_0\oplus\E=(\Sym_0^{(1)}\oplus\Sym_0^{(2)}\oplus\dots)\oplus\E$$
after finitely many steps.\end{proof}

Note this Lemma fails when $x_0=0$, since then $\Sym_0=\Sym$.

\begin{lemma}
\label{diffeo}
There is a neighborhood $U$ of $x_0$ in $\Sym$ and
an analytic map $X:U\to\g_0^\perp$ such that $e^{X(x)}xe^{-X(x)}$ lies in $\Sym_0$ for $x\in U$.
\end{lemma}

\begin{proof}
For $x\in \Sym\oplus\g$, let $\ad(x):\Sym\oplus\g\to\Sym\oplus\g$ be bracketing with $x$, $\ad(x)(y)=[x,y]$.
Then $\ad(x)$ preserves the decomposition $\Sym\oplus\g$ if $x\in\g$ and reverses it if $x\in\Sym$.
If $\langle x,y\rangle=\Trace(xy^t)$ is the usual inner product on $\Sym\oplus\g$, then 
$$\langle [x,y],z\rangle=\langle x,[z,y]\rangle$$
when $x,y\in\Sym$ and $z\in\g$.  This implies that the adjoint of $\ad(x):\Sym\to\g$ is 
$-\ad(x):\g\to\Sym$ when $x\in\Sym$. 

Let $\g_0^\perp$ denote the orthogonal complement of $\g_0$ in $\g$ and let $\Sym_0^\perp$ 
denote the orthogonal complement of 
$\Sym_0$ in $\Sym$. Since the null-space of $\ad(x_0):\Sym\to\g$ equals $\Sym_0$, it follows that the range 
of $\ad(x_0):\g\to\Sym$ is $\Sym_0^\perp$. Since $[\g_0,x_0]=0$, we conclude $[\g_0^\perp,x_0]=\Sym_0^\perp$.

Define a map $\g_0^\perp\oplus\Sym_0\to \Sym$ by 
$$(X,x)\mapsto e^{-X}xe^{X}.$$
At $(0,x_0)$, the derivative of this map is the linear map $(X,x)\mapsto (-[X,x_0])\oplus x$, whose range
equals $\Sym_0^\perp\oplus\Sym_0=\Sym$.
Thus the map is a diffeomorphism at $(0,x_0)$ onto a neighborhood $U$ of $x_0$ in $\Sym$; inverting this
map, the result follows.\end{proof}

\begin{lemma}
\label{S0+diffeo}
$F$ is $C^{n,N}$ on $\Sym^*\times\E$.
\end{lemma}

\begin{proof}
Combining  the two previous lemmas shows $F$ is $C^{n,N}$ on $U\times\E$ hence on $\Sym^*\times\E$.
\end{proof}

At this point that we are left with establishing smoothness near $x_0=0$; this case is more 
significant than at first appears 
as the proof of Lemma \ref{S0} shows that the zero matrix ``propagates'' into larger
and larger subspaces of $\Sym$. Nevertheless, we may be more specific about the asymptotic 
behavior of $F|_{\D\oplus\E}$ at the zero matrix:

\begin{lemma}
\label{WLOG}
Without loss of generality, we may assume in addition that $D^k(F|_{\D\oplus\E})=o(|x|^{n-k})$ as 
$x\to0$ in $\D$ for $0\le k\le n$.
\end{lemma}

\begin{proof}Let $t$ be the $n$-th order Taylor polynomial  of $F|_{\D\oplus\E}$ centered at $x_0=0$.
Since $F$ is rotation-invariant, $F|_{\D\oplus\E}$ is permutation-invariant, hence \cite{MR2003e:00003} 
there is a $C^{\infty,N}$
function $p$ on $\D\oplus\E$, polynomial on $\D$, such that $t=p\circ n$, where $n=(n_1,\dots,n_d)$ are
the Newton sums. Since the Newtons sums extend to polynomial functions on $\Sym\oplus\E$, $t$ extends to
a polynomial function $T$ on $\Sym\oplus\E$; replacing $F$ by $F-T$, we are done.
\end{proof}

To establish smoothness at the origin, we derive a representation formula for $D^nF$ in terms of derivatives
of $F|_{\D\oplus\E}$, which is also of independent interest. This representation formula involves passing
from the coordinate $x\in\Sym$ to ``polar coordinates'' $(r,\pi)$ with $r\in\D$ in a manner analogous to
that presented in the previous section.

A {\it projection} is a real $d\times d$ symmetric matrix $\pi$ satisfying $\pi^2=\pi$, 
and a {\it flag} is a $d$-tuple  $\pi=(\pi_1,\dots,\pi_d)$ of {\it one-dimensional} 
projections that are mutually orthogonal, 
$\pi_i\pi_j=0$ for $i\not=j$, and sum to the identity $\sum_i\pi_i=I$. Since $\Trace(\pi_i)=1$, 
$\pi_i$ is not in $\Sym$.  

Given $x\in\Sym$, let $r_1,\dots,r_d$ denote  its eigenvalues, listed with 
multiplicity, and let $\pi_1,\dots,\pi_d$ denote the projections onto a corresponding 
orthonormal basis of eigenvectors.
Then $r=(r_1,\dots,r_d)$ is in the space $\R^d_0$ of vectors satisfying $r_1+\dots+r_d=0$
and $\pi=(\pi_1,\dots,\pi_d)$ is a flag.  Conversely, if $r\in\R^d_0$ and $\pi$ is a flag,
\beq
x=\sum_i r_i\pi_i
\eeq{decomp}
is in $\Sym$. It is easy to see that the set $\F\subset \Sym^d$ of flags is a compact metric space.

We say a flag $\pi=(\pi_1,\dots,\pi_d)$
is an {\it eigenflag} of $x$ if \rf{decomp} holds for some vector $r$; this happens iff 
$x\pi_i=\pi_i x=r_i\pi_i$ for $i=1,\dots,d$. 

Let $(\R^d)'$ denote the open dense subset of vectors in $\R^d$ with distinct entries and let
$\Sym'$ denote the subset of traceless symmetric matrices with distinct eigenvalues. 

Let $r=r(x)$ equal to the vector of eigenvalues of $x\in \Sym$, arranged in decreasing order;
using the compactness of $\F$, it follows easily that $r:\Sym\to\R^d_0$ is continuous. If $x\in \Sym'$, the 
corresponding eigenflag $\pi=\pi(x)$ is uniquely determined; this is not so if $x$ has repeated eigenvalues. 
We claim the maps $x\mapsto r_i(x)$, $x\mapsto\pi_i(x)$  are analytic on $\Sym'$, 
and we compute the derivatives 
$r_{i\xi}$ and $\pi_{i\xi}$, $i=1,\dots,d$, in the direction of $\xi\in\Sym$; 
this is a standard computation \cite{MR34:3324}.

Let $n_k:\R^d\to\R$ be the $k$-th newton sum, $n_k(r)=(r_1^k+\dots+r_d^k)/k$,
and let $n:\R^d\to\R^d$ be $n=(n_1,\dots,n_d)$. Also define $n:\Sym\to\R^d$ by $n=(n_1,\dots,n_d)$ with 
$n_k(x)=\Trace(x^k)/k$. Then $n(x)=n(r(x))$. Since 
$$\det(Dn(r))=\det
\begin{pmatrix} 1&1&\dots&1\\
r_1&r_2&\dots& r_d\\
&&\ddots&\\
r_1^{d-1}&r_2^{d-1}&\dots&r_d^{d-1}
\end{pmatrix}$$ 
is the Vandermonde determinant, 
$n:(\R^d)'\to\R^d$ is a local diffeomorphism. Since 
$r=n^{-1}(n(r))=n^{-1}(n(x))$, we conclude $r$ is analytic on $\Sym'$.

\begin{lemma}
For $\xi\in\Sym$, we have
\beq
r_{i\xi}=\langle \pi_i,\xi\rangle,\qquad
\pi_{i\xi}=\sum_{j\not=i}\frac{\pi_j\xi \pi_i+\pi_i\xi \pi_j}{r_i-r_j},
\eeq{pdiff}
$i=1,\dots,d$, on $\Sym'$.
\end{lemma}

\begin{proof}
By the chain rule, 
$$Dr(x)=D (n^{-1}\circ n)(x)=D(n^{-1})(n(r)) \cdot Dn(x)=(Dn(r))^{-1}\cdot Dn(x).$$
Since $Dn_k(x)=x^{k-1}$, Cramer's rule yields
$$Dr(x)=\begin{pmatrix} \pi_1\\ \pi_2\\ \dots\\ \pi_d\end{pmatrix},$$
or, what is the same, 
\beq
r_{i\xi}=\langle\pi_i,\xi\rangle,\qquad i=1,\dots,d.
\eeq{rdiff}
In particular, since $r_i$ is analytic, this shows that 
the maps $\pi_i$, $i=1,\dots,d$, are  analytic.

To compute $\pi_{i\xi}$, differentiate $x\pi_i=\pi_ix=r_i\pi_i$ to get
$\xi \pi_i+x\pi_{i\xi}=r_{i\xi}\pi_i+r_i\pi_{i\xi}$.
Left multiply by $\pi_j$, $j\not=i$, to get
$$\pi_j\xi \pi_i+r_j\pi_j\pi_{i\xi}=
\pi_j\xi \pi_i+\pi_j x \pi_{i\xi}=\pi_jr_{i\xi}\pi_i+\pi_jr_i\pi_{i\xi}=r_i\pi_j\pi_{i\xi}$$
which yields
$$\pi_j\pi_{i\xi}=\frac{\pi_j\xi \pi_i}{r_i-r_j}=\pi_j\pi_{i\xi}\pi_i.$$
Differentiating $\pi_i^2=\pi_i$, we obtain $\pi_i\pi_{i\xi}+\pi_{i\xi}\pi_i=\pi_{i\xi}$, hence
$\pi_i\pi_{i\xi}\pi_i=0$; summing over $j$, we conclude
$$\pi_{i\xi}\pi_i=\sum_{j\not=i}\frac{\pi_j\xi \pi_i}{r_i-r_j}.$$
Adding this last equation to its transpose, we arrive at  \rf{pdiff}.
\end{proof}

We say a  function $f:\R^d_0\times\F\times\E\to\R$ is {\it symmetric} if
$$
f(r_{\sigma1},\dots,r_{\sigma d},\pi_{\sigma1},\dots,\pi_{\sigma d},v)=
f(r_1,\dots,r_d,,\pi_1,\dots,\pi_d,v)$$ 
holds on $\R^d_0\times\F\times\E$ for every permutation $\sigma$. A subset $K$ is {\it symmetric}
if ${\bf 1}_K$ is symmetric.

We say $f$ is {\it consistent} if $f(r,\pi,v)=f(r,\pi',v)$ whenever $x=\sum_ir_i\pi_i=\sum_ir_i\pi'_i$.
This is the same as saying $f(r,\pi,v)=f(r,\pi',v)$ whenever $\pi_\lambda=\pi_\lambda'$, where
$$\pi_\lambda=\sum_{r_i=\lambda} \pi_i$$
for every eigenvalue $\lambda$ of $x$. 

If $F:\Sym\times\E\to\R$ is a continuous function, then $f:\R^d_0\times\F\times\E\to\R$ defined by
\beq
f(r,\pi,v)=F\left(\sum_i r_i\pi_i,v\right),\qquad r\in\R^d_0,\pi\in \F,v\in\E,
\eeq{char}
is clearly continuous, symmetric and consistent. Note that $F$ is rotation-invariant iff $f$ does not
depend on $\pi$.

\begin{lemma}
\label{consist}
If $f:\R^d_0\times\F\times\E\to\R$ is continuous, symmetric, and consistent,
there exists a unique continuous $F:\Sym\times\E\to\R$ satisfying \rf{char}.
\end{lemma}

\begin{proof}
Since $r_i$, $\pi_i$, $i=1,\dots,d$, are analytic on $\Sym'$ and $f$ is symmetric,
it is clear that \rf{char} defines $F$ uniquely and continuously on $\Sym'\times\E$.
If $x_n\in\Sym'$ and $x_n\to x\in\Sym$ and $v_n\to v$ in $\E$, 
we need to establish the convergence of $(F(x_n,v_n))$.
To this end, let $r_n$ denote the corresponding vectors of eigenvalues, arranged in decreasing
order, and let $\pi_n$ denote the corresponding eigenflags.  Then $r_n$ converges to the vector $r$
of eigenvalues of $x$ arranged in non-increasing order. 
If $\pi$ is a limit point of $(\pi_n)$, then $\pi$ is an eigenflag of $x$.
By consistency, $f(r,\pi,v)$ depends only on $r$ and the projections $\pi_\lambda$
onto the $\lambda$-eigenspaces of $x$, hence only on $x$. Thus $f(r,\pi,v)$ 
does not depend on the subsequence,
$(f(r_n,\pi_n,v_n))=(F(x_n,v_n))$ converges to a limit, and \rf{char} holds at all $(r,\pi,v)$.
\end{proof}

Let $\Sym^d=\Sym\times \Sym\dots\times\Sym$ be the $d$-fold product.
Given $\xi\in\Sym$ and a skew-symmetric $d\times d$ matrix $a$, define a map 
$\delta=\delta(a,\xi):\Sym^d\to\Sym^d$ by
$$\delta(a,\xi)(\pi)_i=\pi_i\xi\left(\sum_j a_{ij}\pi_j\right)+
\left(\sum_j a_{ij}\pi_j\right)\xi\pi_i,\qquad i=1,\dots,d.$$

\begin{lemma}
The map $\delta$ restricted to $\F$ is a vector field tangent to $\F$.
\end{lemma}

\begin{proof}
To see this, let $\pi(t)\in\Sym^d$ be a smooth curve of $d$-tuples of symmetric matrices 
starting at $\pi(0)\in\F$ satisfying 
$\dot\pi_i=\delta_i(\pi)$, $i=1,\dots,d$, for $t$ small. We show $\pi(t)\in\F$ by showing
\begin{enumerate}
\item $\sum_i\pi_i(t)=1$,
\item $\pi_i(t)\pi_j(t)=0$, $i\not=j$,
\item $\pi_i(t)^2=\pi_i(t)$.
\end{enumerate}
(1) follows since $\sum_i\dot\pi_i(t)=\sum_i\delta_i(\pi(t))=0$.
Differentiation shows that $x_{ij}(t)=\pi_i(t)\pi_j(t)$, $i\not=j$, satisfies a linear system of 
differential equations with time-varying coefficients; since $x_{ij}(0)=0$, $i\not=j$, (2) follows.
(3) follows since by (2) $(d/dt)\pi_i(t)^2=\dot\pi_i(t)\pi_i(t)+\pi_i(t)\dot\pi_i(t)=\dot\pi_i(t)$.
Thus $\pi(t)\in\F$. 
\end{proof}

Define vector fields $\delta_{ij\xi}=-\delta_{ji\xi}$, $i\not=j$, on $\F$ by
$$\delta(a,\xi)=\sum_{i\not=j}a_{ij}\delta_{ij\xi}.$$
Then for each $i,j,\xi$, $\delta_{ij\xi}$ is a vector field on $\F$.
Note that \rf{pdiff} can be rewritten as
$$\pi_\xi=\frac12\sum_{i\not=j}\frac{\delta_{ij\xi}(\pi)}{r_i-r_j}=\frac12\delta
\left(\rho,\xi\right),$$
where $\rho$ is the skew-symmetric matrix with entries $1/(r_i-r_j)$.

Let $(\R^d_0)'$ denote the vectors in $\R^d_0$ with distinct entries and
let $(\R^d_0)^*$ be the nonzero vectors in $\R^d_0$. 

If $f:\R^d_0\times\F\times\E\to\R^d$ is polynomial in $\pi$, $f=(f_1,\dots,f_d)$, let
$$\L_\xi(f)(r,\pi,v)=\sum_if_{i}\cdot\langle\pi_i,\xi\rangle+\frac14\sum_{i\not=j}
\int_0^1 \delta_{ij\xi}(f_{i}-f_{j})(r(t),\pi,v)\,dt.$$
If $f$ is $C^{n,N}$ on $\R^d_0\oplus\E$ and polynomial in $\pi$, so is $\L_\xi(f)$.

Let $Df=(f_{r_1},\dots,f_{r_d})$ be the gradient of $f$ in $r$.

\begin{lemma}
If $f$ given by \rf{char} is $C^{n,N}$ on $\R^d_0\oplus\E$ and polynomial in $\pi$,
then
\beq
D_\xi^n F\left(\sum_ir_i\pi_i,v\right)=(\L_\xi D)^nf(r,\pi,v)
\eeq{der}
on $(\R^d_0)'\times\F\times\E$.
\end{lemma} 
\begin{proof}
Let $\pi(t)$ be the integral curve of $\delta_{ij\xi}$ starting from $\pi\in\F$; since the sum of
$i$-th and $j$-th components of $\delta_{ij\xi}$ vanishes, $\pi_i(t)+\pi_j(t)$ does not depend on $t$;
then $f$ consistent and polynomial in $\pi$
and $r_i=r_j$ implies $f(r,\pi(t),v)$ does not depend on $t$, hence $\delta_{ij\xi}(f)(r,\pi)=0$.

Given $r\in\R^d$, let $r(t)$ differ from $r$ only in the $i$-th and $j$-th components, by setting
$r_i(t)=tr_i+(1-t)(r_i+r_j)/2$, $r_j(t)=tr_j+(1-t)(r_i+r_j)/2$.

Since $\delta_{ij\xi}(f)(r(0),\pi,v)$ vanishes,
the fundamental theorem of calculus applied to $\delta_{ij\xi}(f)(r(t),\pi,v)$ implies
\beq
\frac{\delta_{ij\xi}(f)(r,\pi,v)}{r_i-r_j}=\frac12\int_0^1 \delta_{ij\xi}(f_{r_i}-f_{r_j})(r(t),\pi,v)\,dt.
\eeq{int}
By the chain rule,\rf{int}, and \rf{pdiff},  
\beq
\begin{split}
D_\xi F\left(\sum_ir_i\pi_i,v\right)&=\sum_if_{r_i}\cdot\langle\pi_i,\xi\rangle+\frac12\sum_{i\not=j}
\frac{\delta_{ij\xi}(f)}{r_i-r_j}\\
&=\L_\xi(D(f))(r,\pi,v)
\end{split}
\eeq{LL}
on $(\R^d_0)'\times\F\times\E$. If $F$ is $C^{1,N}$, \rf{LL} is valid on $\R^d_0\times\F\times\E$, thus 
$$D^2_\xi F\left(\sum_ir_i\pi_i,v\right)=(\L_\xi D)^2(f)(r,\pi,v)$$
 on $(\R^d_0)'\times\F\times\E$. 
Since $F$ is $C^{n-1,N}$ we may repeat this argument $n-1$ times;
the result follows.
\end{proof}

If $F$ is rotation-invariant, then $F$ is $C^{n,N}$ on $\Sym^*\oplus\E$, and hence \rf{der} is valid on
$(\R^d_0)^*\times\F\times\E$.  Moreover, $f=f(r,v)$ does not depend on $\pi$ and hence \rf{der} implies
\beq
\left|D^n F(x,v)\right|\le C\sup_{|r|\le|x|}|D^nf(r,v)|.
\eeq{decay}
Recalling Lemma \ref{WLOG}, this implies $D^n F(x,v)\to0$ as 
$|x|=|r|\to0$; since we know $F$ is $C^{{n-1},N}$, this implies $F(\cdot,v)$ is $C^{n}$ on $\Sym$ for each $v\in\E$.
This in turn implies the validity of \rf{der}  on $\R^d_0\times\F\times\E$, which in turn
implies $F$ is $C^{n,N}$ on $\Sym\oplus\E$.  This completes the proof of Theorem \ref{main}.

We now prove the Corollary.
The first statement is an immediate consequence of Theorem \ref{diffth}.
Away from the origin, if $f$ is $C^{n,\alpha}$ or analytic, the proof (Lemmas \ref{S0}, \ref{diffeo},
\ref{S0+diffeo}) of
Theorem \ref{main}, unchanged, establishes $F$ is $C^{n,\alpha}$ or analytic respectively.
If $f$ is $C^{n,\alpha}$, then, from Theorem \ref{diffth}, $F$ is $C^n$. If $t_n$ is the $n$-th order
Taylor polynomial of $f$ at the origin, then $t_n$ is permutation-invariant, hence $t_n$ is the
restriction to $\D$ of a rotation-invariant polynomial $T_n$ on $\Sym$. It follows that $T_n$ is the $n$-th
order Taylor polynomial of $F$ at the origin. Replacing $F$ by $F-T_n$, since $D^nt_n(r)=D^nf(0)$ and 
$D^nT_n(x)=D^nF(0)$,  by \rf{decay} we have
$$\left|D^n F(x)-D^nF(0)\right|\le C\sup_{|r|\le|x|}|D^nf(r)-D^nf(0))|.$$
Thus $F$ is $C^{n,\alpha}$ at the origin.
If $f$ is analytic at the origin, $|f(r)-t_n(r)|\le C|r/2\epsilon|^n$ on $|r|<\epsilon$ 
for all $n$; since $F$, $T_n$ and $|x|^n$ are rotation-invariant,
it follows that $|F(x)-T_n(x)|\le C|x/2\epsilon|^n$ on $|x|<\epsilon$ for all $n$; 
thus $F$ is analytic at the origin.

\end{section}

\bibliographystyle{amsplain}

\begin{thebibliography}{99999}
\bibitem{MR19:832b}
Chandler Davis.
\newblock All convex invariant functions of hermitian matrices.
\newblock {\em Arch. Math.}, 8:276--278, 1957.

\bibitem{GH}
Yury Grabovsky and Omar Hijab.
\newblock A generalization of the Chandler Davis theorem.
\newblock 2003, submitted.

\bibitem{MR34:3324}
Tosio Kato.
\newblock {\em Perturbation theory for linear operators}.
\newblock Die Grundlehren der mathematischen Wissenschaften, Band 132.
  Springer-Verlag New York, Inc., New York, 1966.

\bibitem{MR2003c:22001}
Anthony~W. Knapp.
\newblock {\em Lie groups beyond an introduction}, volume 140 of {\em Progress
  in Mathematics}.
\newblock Birkh\"auser Boston Inc., Boston, MA, second edition, 2002.

\bibitem{MR2003e:00003}
Serge Lang.
\newblock {\em Algebra}, volume 211 of {\em Graduate Texts in Mathematics}.
\newblock Springer-Verlag, New York, third edition, 2002.

\bibitem{MR2001g:49026}
A.~S. Lewis.
\newblock Convex analysis on {C}artan subspaces.
\newblock {\em Nonlinear Anal.}, 42(5, Ser. A: Theory Methods):813--820, 2000.

\bibitem{MR2002i:15013}
Adrian~S. Lewis and Hristo~S. Sendov.
\newblock Quadratic expansions of spectral functions.
\newblock {\em Linear Algebra Appl.}, 340:97--121, 2002.

\bibitem{math.FA/0208223}
Igor Rivin.
\newblock {Another Simple Proof of a Theorem of Chandler Davis,\/}
\newblock arxiv.org preprint math.FA/0208223
\newblock Submitted.

\end{thebibliography}

\end{document}